\documentclass[11pt,amsfonts]{article}
\usepackage{graphicx}
\usepackage{latexsym}
\usepackage{amssymb}
\usepackage{amsmath}
\usepackage{enumerate}
\usepackage{stmaryrd}
\usepackage{layout}
\usepackage{eufrak}
\newtheorem{prop}{Proposition}

\newtheorem{theorem}{Theorem}
\newtheorem{remark}{Remark}

\def\real{{\mathord{{\rm I\kern-2.8pt R}}}}        
\def\inte{{\mathord{{\rm I\kern-2.8pt N}}}}

\def\sZZ{{\rm Z\kern-2.8ptem{}Z}}

\def\z{{\mathchoice
  {\sZZ}
  {\sZZ}
  {\rm Z\kern-0.30em{}Z}
  {\rm Z\kern-0.25em{}Z} }}
\def\sQQ{{\kern 0.27em \vrule height1.45ex width0.03em depth0em
          \kern-0.30em \rm Q}}
\def\qu{{\mathchoice
    {\sQQ}
    {\sQQ}
  {\kern 0.225em \vrule height1.05ex width0.025em depth0em \kern-0.25em \rm Q}
  {\kern 0.180em \vrule height0.78ex width0.020em depth0em \kern-0.20em \rm Q}
        }}
\def\sCC{{\kern 0.27em \vrule height1.45ex width0.03em depth0em
          \kern-0.30em \rm C}}
\def\complex{{\mathchoice
    {\sCC}
    {\sCC}
  {\kern 0.225em \vrule height1.05ex width0.025em depth0em \kern-0.25em \rm C}
  {\kern 0.180em \vrule height0.78ex width0.020em depth0em \kern-0.20em \rm C}
        }}

\newcommand{\ba}{\begin{array}}
\newcommand{\ea}{\end{array}}
\newcommand{\be}{\begin{equation}}
\newcommand{\ee}{\end{equation}}
\newcommand{\bea}{\begin{eqnarray}}
\newcommand{\eea}{\end{eqnarray}}
\newcommand{\beaa}{\begin{eqnarray*}}
\newcommand{\eeaa}{\end{eqnarray*}}

\def\z{\zeta}

\font\tenmath=msbm10 \font\sevenmath=msbm7 \font\fivemath=msbm5
\newfam\mathfam \textfont\mathfam=\tenmath
\scriptfont\mathfam=\sevenmath \scriptscriptfont\mathfam=\fivemath

\def \={{\buildrel {\rm (law)} \over =}}

\def\qed{ \hfill \vrule width.25cm height.25cm depth0cm\smallskip}

\newcommand{\basa}{\begin{assumption}}
\newcommand{\easa}{\end{assumption}}

\newcommand{\bas}{\begin{assum}}
\newcommand{\eas}{\end{assum}}

\newcommand{\ignore}[1]{}
\begin{document}

\renewcommand{\thefootnote}{\fnsymbol{footnote}}

\renewcommand{\thefootnote}{\fnsymbol{footnote}}
\begin{center}
{\large{\textsc{On the law of the solution to a stochastic heat equation with fractional noise in time \footnote{\noindent Dedicated to the memory of Constantin Tudor}}}}\\~\\
Solesne Bourguin\footnote{\noindent SAMM, Universit\'e de Paris 1 Panth\'eon-Sorbonne, 90, rue de Tolbiac, 75634, Paris, France.
Email: {\tt solesne.bourguin@univ-paris1.fr}
}
and
Ciprian A. Tudor\footnote{\noindent Laboratoire Paul Painlev\'e, Universit\'e de Lille 1, F-59655 Villeneuve d'Ascq, France.
Email: {\tt tudor@math.univ-lille1.fr}}\footnote{\noindent Associate member of the team Samm, Universit\'e de Paris 1 Panth\'eon-Sorbonne. The second  author is partially supported by the ANR grant "Masterie" BLAN 012103.}\\
{\it Universit\'e Paris 1} and {\it Universit\'e Lille 1}\\~\\
\end{center}
{\small \noindent {\bf Abstract:}
We study the law of the solution to the stochastic heat equation with additive Gaussian noise which behaves as the fractional Brownian motion in time and is white in space. We prove a decomposition of the solution  in terms of  the bifractional Brownian motion. Our result is an extension of a result in the paper \cite{Swanson1}.\\


\noindent
{\bf 2010 AMS Classification Numbers:} 60H15, 60H05. }
\\
\section{Introduction}
The purpose of this note is to make some remarks on the law of the solution to a heat equation with a fractional noise in time. Let us briefly describe the context and the motivation.  Consider a centered Gaussian field $\left( W(t,x), t\in [0,T] , x\in \mathbb{R} ^{d}\right)$ with covariance given by
\begin{equation}
\label{cov1}
\mathbf{E}\left(  W(t,x)W(s,y) \right) = (s\wedge t) (x\wedge y)
\end{equation}
where for $x=(x_{1},..,x_{d}), y= (y_{1},..., y_{d}) \in \mathbb{R} ^{d}$ and where we denoted $x\wedge y= \prod_{i=1}^{d} \left(x_{i} \wedge y_{i}\right)$.  Also consider the stochastic partial differential equation with additive noise
\begin{eqnarray}
\label{eq1} \frac{\partial u}{\partial t} &=& \frac{1}{2}
\Delta u + \dot {W} , \quad t  \in [0,T], \mbox{\ } x\in \mathbb{R}^{d} \\
\nonumber u_{0,x} &=& 0, \quad x \in \mathbb{R}^d,
\end{eqnarray}
where the noise $W$ is defined by (\ref{cov1}). The noise $W$ is usually refered to as a {\it space-time white noise}. It is well-known (see for example the seminal paper by Dalang \cite{Da}) that the heat equation (\ref{eq1}) admits a unique mild solution if and only if $d=1$. This mild solution is defined as
\begin{equation}
\label{sol}
 u(t,x)=\int_{0}^{T}
\int_{\mathbb{R}^{d}} G(t-s,x-y)W(ds,dy), \hskip0.5cm t\in [0,T], x\in \mathbb{R} ^{d}
\end{equation}
where the above integral is a Wiener integral with respect to the Gaussian process $W$ (see e.g. \cite{BalanTudor1} for details) and $G $ is the Green kernel of the heat equation given by
\begin{equation}
\label{fund-sol-heat} G(t,x)=\left\{
\begin{array}{ll} (4 \pi t)^{-d/2} \exp\left( -\frac{|x|^{2}}{4t}\right) & \mbox{if $t>0, x \in \mathbb{R}^{d}$} \\
0 & \mbox{if $t \leq 0, x \in \mathbb{R}^{d}$}.
\end{array} \right.
\end{equation}
 Consequently, the process $(u(t,x), t\in [0,T] , x \in \mathbb{R}) $ is a centered Gaussian process. The question we address in this paper comes from the following observation: it has been proven in \cite{Swanson1} for $d=1$ that the covariance of the solution (\ref{sol}) satisfies the following assertion :  for every $x\in \mathbb{R}$ we have
\begin{equation}\label{covu}
\mathbf{E} \left(  u(t,x) u(s,x) \right) = \frac{1}{\sqrt{2\pi} } \left( \sqrt{ t+s} -\sqrt{\vert t-s\vert } \right), \mbox{ for every } s,t \in [0,T].
\end{equation}
This fact establishes an interesting connection between the law of the solution (\ref{sol}) and the so-called {\it bifractional Brownian motion.}   Recall that the bifractional Brownian motion $(B^{H,K}_{t})_{t \in [0,T]}$ is a centered Gaussian process, starting from zero, with covariance
  \begin{equation}
    \label{cov-bi}
    R^{H,K}(t,s) := R(t,s)= \frac{1}{2^{K}}\left( \left(
        t^{2H}+s^{2H}\right)^{K} -\vert t-s \vert ^{2HK}\right), \hskip0.5cm s,t \in [0,T]
  \end{equation}
  with $H\in (0,1)$ and $K\in (0,1]$. We refer to \cite{HV} and \cite{RT} for the definition and the basic properties of this process. Note that, if $K=1$ then $B^{H,1}$ is a fractional Brownian motion with Hurst parameter $H\in
  (0,1)$. When $K=1$ and $H=\frac{1}{2}$ then it reduces to the standard Brownian motion.  Relation (\ref{covu}) implies that the solution $\left( u(t,x), t \in [0,T], x \in \mathbb{R}\right)  $ to equation (\ref{eq1}) is a bifractional Brownian motion with parameters $H=K=\frac{1}{2}$ multiplied by the constant $2^{-K} \frac{1}{\sqrt{2\pi }}$. This is not the only connection between the solution of the heat equation and fractional processes. Another interesting link has been proven in \cite{M}.  

\noindent Our purpose is to study the law of the linear heat equation driven by a fractional noise in time. That is, we will consider a Gaussian field
$\left( W^{H} (t,x), t\in [0,T] , x\in \mathbb{R} ^{d} \right) $ with covariance
\begin{equation}
\label{cov2}
\mathbf{E} \left(  W^{H} (t,x) W^{H} (s, y) \right) = \frac{1}{2} (t^{2H}+ s^{2H} -\vert t-s \vert ^{2H} ) (x\wedge y)
\end{equation}
where we will assume throughout the paper that the Hurst parameter $H$ is contained in the interval $(\frac{1}{2}, 1)$.
\noindent Let us consider the linear stochastic heat equation
\begin{eqnarray}
\label{heatEq}
u_{t} = \Delta u + \dot{W^{H}}, \mbox{\  \  \  } t \in \left[0,T\right], x \in \mathbb{R}^{d}
\end{eqnarray}
with $u(.,0) = 0$, where $(W(t,x))_{t \in \left[0,T\right], x \in \mathbb{R}^{d}}$ is a centered Gaussian noise with covariance (\ref{cov2}). The solution of (\ref{heatEq}) can be written in {\it mild form}  as
\begin{eqnarray}\label{n100}
U(t,x) = \int_{0}^{t}\int_{\mathbb{R}^{d}}G(t-s,x-y)W^{H}(ds,dy)
\end{eqnarray}
where the above integral is a Wiener integral with respect to the noise $W^{H}$ (see e.g \cite{BalanTudor1} for the definition) and $G$ is given (\ref{fund-sol-heat}).

\noindent The following result has been proven in \cite{BalanTudor1} :

\begin{theorem}
The process $(U (t,x))_{t \in \left[0,T\right], x \in \mathbb{R}^{d}}$ exists and satisfies $$ \underset{t \in \left[0,T\right], x \in \mathbb{R}^{d}}{\mbox{sup}}\mathbf{E}\left(U(t,x)^{2}\right) < +\infty $$ if and only if $d<4H$.
\end{theorem}

\noindent This implies that, in contrast to the white-noise case, we are allowed to consider the spatial dimension $d$ to be $1,2$ or 3.

\noindent We are interested in characterizing the law of the Gaussian process $U(t,.)$ defined by (\ref{n100}) and to see its relation with the bifractional Brownian motion.  In other words, we try to understand how relation (\ref{covu}) changes when the noise becomes fractional in time. We prove that the solution to the heat equation with additive noise fractional in time and white in space can be decomposed in terms of the bifractional Brownian motion with Hurst parameters $H=K=\frac{1}{2}$ but in this case some  additional Gaussian processes appear. Concretely, we show that  for every $x\in \mathbb{R} ^{d}$,
\begin{equation*}
\left( U(t,x) + Y^{H}, t\in [0,T]\right)  \stackrel{\rm Law}{=}\left(  C_{0} B^{\frac{1}{2}, \frac{1}{2}}_{t} + X^{H}_{t}, t\in [0,T]\right),
\end{equation*}where $C_{0}>0$ is an explicit constant and $X^{H}, Y^{H}$ are independent Gaussian processes, also independent from $U$, with covariance functions (\ref{rx}) and (\ref{ry}) respectively. When $H=\frac{1}{2}$, an interesting phenomenon happens: the process    $Y^{H}$ vanishes while the constant $C_{0}$ converges to $\frac{1}{2 \sqrt{2\pi}}$ when $H$ is close to one half and the  process $X^{H}$ converges to $\frac{1}{2 \sqrt{2\pi}}B^{\frac{1}{2}, \frac{1}{2}}$ for $H$ close to one half (see Remark \ref{onehalf}). In other words, in the fractional case $H\not= \frac{1}{2}$ the solution ``keeps'' half of the bifractional Brownian motion $B^{\frac{1}{2}, \frac{1}{2}}$ while the other half ``spreads'' into two parts $X^{H}$ and $Y^{H}$.

\section{Analysis of the covariance of the process $U(t,.)$}

Suppose that $s \leq t$ and denote by
\begin{eqnarray*}
R(t,s) = \mathbf{E}\left( U(t,x)U(s,x) \right)
\end{eqnarray*}
where $x \in \mathbb{R}^{d}$ is fixed. We can write, from the proof of Theorem 2.7 in \cite{BalanTudor1},
\begin{eqnarray}
\label{covar}
R(t,s) = \alpha_{H}\frac{1}{2\sqrt{2\pi}}\int_{0}^{t}\int_{0}^{s}\left|u-v\right|^{2H-2}((t+s) - (u+v))^{-\frac{d}{2}}dvdu
\end{eqnarray}
with $\alpha_{H} = H(2H-1)  $.


\noindent The purpose of this paragraph is to analyze the covariance (\ref{covar}) of the solution $U(t,x)$ and to understand its relation with the bifractional Brownian motion.
The following proposition gives a decomposition of the covariance function of $U(t,.)$ in the case $d\neq 2$. The lines of the proof below will explain why the case $d=2$ has to be excluded.
\begin{prop}\label{pp1}
Suppose $d\not=2$. The covariance function $R(t,s)$ can be decomposed as follows
\begin{eqnarray*}
R(t,s) = \frac{1}{2 \sqrt{2\pi}}\alpha_{H}C_{d}\beta\left(2H-1,-\frac{d}{2} + 2\right)\left[(t+s)^{2H-\frac{d}{2}} - (t-s)^{2H-\frac{d}{2}}\right] + R_{1} ^{(d)}(t,s)
\end{eqnarray*}
where $C_{d} = \frac{2}{2-d}$, $\beta(x,y)$ is the Beta function defined for $x,y >0$ by $\beta(x,y) = \int_{0}^{1}t^{x-1}(1-t)^{y-1}dt$ and
\begin{eqnarray*}
R_{1} ^{(d)}(t,s) &=& \frac{1}{2 \sqrt{2\pi}}\alpha_{H}C_{d}\left[\int_{0}^{s}daa^{2H-2}\left[((t+s) - a)^{-\frac{d}{2}+1} - ((t-s) + a)^{-\frac{d}{2}+1}\right]  \right. \nonumber \\
&&\left. -\int_{0} ^{s} da (s-a) ^{-\frac{d}{2}+1} \left[ (t-a)^{2H-2} + (t+a) ^{2H-2}\right]\right].
\end{eqnarray*}
\end{prop}
{\noindent \bf Proof: }Fix $t>s$. By performing the change of variables $u-v = a$ and $u+v = b$ with $a+b = 2u \in (0,2t)$ and $b-a = 2v \in (0,2s)$ in (\ref{covar}), we get
\begin{eqnarray*}
R(t,s) &=& \frac{1}{2 \sqrt{2\pi}}\alpha_{H}\int_{-s}^{t}\left|a\right|^{2H-2}\int_{a \vee (-a)}^{(2t-a) \wedge (2s+a)}((t+s) - b)^{-\frac{d}{2}}dbda
\\
&=& \frac{1}{2 \sqrt{2\pi}}\alpha_{H}\left[\int_{-s}^{0}(-a)^{2H-2}\int_{-a}^{2s+a}((t+s) - b)^{-\frac{d}{2}}dbda \right.
\\
&& \left. + \int_{0}^{t-s}a^{2H-2}\int_{a}^{2s+a}((t+s) - b)^{-\frac{d}{2}}dbda \right.
\\
&&\left. + \int_{t-s}^{t}a^{2H-2}\int_{a}^{2t-a}((t+s) - b)^{-\frac{d}{2}}dbda\right].
\end{eqnarray*}
By performing the change of variables $a \mapsto (-a)$ in the first summand, we get
\begin{eqnarray*}
R(t,s) &=& \frac{1}{2 \sqrt{2\pi}}\alpha_{H}\left[\int_{0}^{s}a^{2H-2}\int_{a}^{2s-a}((t+s) - b)^{-\frac{d}{2}}dbda \right.
\\
&& \left.+ \int_{0}^{t-s}a^{2H-2}\int_{a}^{2s+a}((t+s) - b)^{-\frac{d}{2}}dbda \right.
\\
&& + \left.\int_{t-s}^{t}a^{2H-2}\int_{a}^{2t-a}((t+s) - b)^{-\frac{d}{2}}dbda\right].
\end{eqnarray*}
\begin{remark}
We can notice why the case $d=2$ must be treated separately in the latter equation. The integral with respect to $db$ involves logarithms and it cannot lead to the covariance of the bifractional Brownian motion.
\end{remark}

\noindent By explicitly computing the inner integrals, we obtain
\begin{eqnarray*}
R(t,s) &=& \frac{1}{2 \sqrt{2\pi}}\alpha_{H}C_{d}\left[\int_{0}^{s}a^{2H-2}\left[-((t+s)-b)^{-\frac{d}{2}+1}\right]_{b=a}^{b=2s-a}da \right.
\\
&& \left.+ \int_{0}^{t-s}a^{2H-2}\left[-((t+s)-b)^{-\frac{d}{2}+1}\right]_{b=a}^{b=2s+a}da\right.
\\
&& \left. + \int_{t-s}^{t}a^{2H-2}\left[-((t+s)-b)^{-\frac{d}{2}+1}\right]_{b=a}^{b=2t-a}da\right]
\\
&=& \frac{1}{2 \sqrt{2\pi}}\alpha_{H}C_{d}\left[\int_{0}^{s}a^{2H-2}((t+s)-a)^{-\frac{d}{2}+1}da - \int_{0}^{s}a^{2H-2}((t-s)+a)^{-\frac{d}{2}+1}da\right]
\\
&& + \alpha_{H}C_{d}\left[\int_{0}^{t-s}a^{2H-2}((t+s)-a)^{-\frac{d}{2}+1}da - \int_{0}^{t-s}a^{2H-2}((t-s)-a)^{-\frac{d}{2}+1}da\right]
\\
&& + \alpha_{H}C_{d}\left[\int_{t-s}^{t}a^{2H-2}((t+s)-a)^{-\frac{d}{2}+1}da - \int_{t-s}^{t}a^{2H-2}(a - (t-s))^{-\frac{d}{2}+1}da\right]
\\
&=&\frac{1}{2 \sqrt{2\pi}} \alpha_{H}C_{d}\left[\int_{0}^{t+s}a^{2H-2}((t+s)-a)^{-\frac{d}{2}+1}da - \int_{0}^{t-s}a^{2H-2}((t-s)-a)^{-\frac{d}{2}+1}da\right]
\\
&& + R_{1} ^{(d)}(t,s)
\end{eqnarray*}
where
\begin{eqnarray}\label{r1}
&&R_{1} ^{(d)}(t,s) \nonumber \\
&=& \frac{1}{2 \sqrt{2\pi}}\alpha_{H}C_{d}\left[\int_{0}^{s}a^{2H-2}((t+s) - a)^{-\frac{d}{2}+1}da \right. \nonumber \\
&& \left. - \int_{0}^{s}a^{2H-2}((t-s) + a)^{-\frac{d}{2}+1}da \right. \nonumber
\\
&& \left. - \int_{t-s}^{t}a^{2H-2}(a-(t-s))^{-\frac{d}{2}+1}da - \int_{t}^{t+s}a^{2H-2}((t+s) - a)^{-\frac{d}{2}+1}da\right].
\end{eqnarray}
At this point, we perform the change of variable $a \mapsto  \frac{a}{t+s}$ and we obtain
\begin{eqnarray*}
\int_{0}^{t+s}a^{2H-2}((t+s)-a)^{-\frac{d}{2}+1}da &=& (t+s)^{2H - \frac{d}{2}}\int_{0}^{1}a^{2H-2}(1-a)^{-\frac{d}{2}+1}da
\\
&=& \beta\left(2H-1,-\frac{d}{2} + 2\right)(t+s)^{2H - \frac{d}{2}}
\end{eqnarray*}
and in the same way, with the change of variable $a \mapsto  \frac{a}{t-s}$, we obtain
\begin{eqnarray*}
\int_{0}^{t-s}a^{2H-2}((t-s)-a)^{-\frac{d}{2}+1}da &=& (t-s)^{2H - \frac{d}{2}}\int_{0}^{1}a^{2H-2}(1-a)^{-\frac{d}{2}+1}da
\\
&=& \beta\left(2H-1,-\frac{d}{2} + 2\right)(t-s)^{2H - \frac{d}{2}}.
\end{eqnarray*}
As a consequence, we obtain
\begin{eqnarray*}
R(t,s) = \alpha_{H}C_{d}\beta\left(2H-1,-\frac{d}{2} + 2\right)\left[(t+s)^{2H-\frac{d}{2}} - (t-s)^{2H-\frac{d}{2}}\right] + R_{1} ^{(d)}(t,s)
\end{eqnarray*}
with $R_{1}^{(d)}$ given by (\ref{r1}). Let us further analyze the function denoted by $R_{1} ^{(d)}(t,s)$. Note that  for every $s,t\in [0,T]$
\begin{equation*}
\left( \frac{1}{2 \sqrt{2\pi}}\right) ^{-1}R_{1} ^{(d)}(t,s)  =A(t,s) + B(t,s)
\end{equation*}
where
\begin{eqnarray*}
A (t,s)=  \alpha_{H}C_{d}\left[\int_{0}^{s}a^{2H-2}((t+s) - a)^{-\frac{d}{2}+1}da - \int_{0}^{s}a^{2H-2}((t-s) + a)^{-\frac{d}{2}+1}da \right]
\end{eqnarray*}
and
\begin{equation*}
B(t,s)=  \alpha_{H}C_{d}\left[- \int_{t-s}^{t}a^{2H-2}(a-(t-s))^{-\frac{d}{2}+1}da - \int_{t}^{t+s}a^{2H-2}((t+s) - a)^{-\frac{d}{2}+1}da\right].
\end{equation*}
By the change of variables $a-t =\tilde{a} $, we can express $B$ as
\begin{eqnarray*}
B(t,s)&= &\alpha _{H} C_{d} \left[ -\int_{-s} ^{0} (a+t) ^{2H-2} (a+s) ^{-\frac{d}{2}+1} da -\int_{0} ^{s} (a+t) ^{2H-2}(s-a)  ^{-\frac{d}{2}+1} da\right]\\
&=& -\alpha _{H} C_{d} \int_{0} ^{s} da (s-a) ^{-\frac{d}{2}+1} \left[ (t-a)^{2H-2} + (t+a) ^{2H-2}\right]
\end{eqnarray*}
and the desired conclusion is obtained.\qed

\vskip0.3cm

\noindent Let us point out that the constant $C_{d}$ is positive for $d=1$ and negative for $d=3$. This partially explains why we obtain different decompositions in these two cases. Thanks to the decomposition in Proposition \ref{pp1}, we can prove the following proposition.

\begin{theorem}\label{tt2} Assume $d=1$ and let $U$ be the solution of the heat equation with fractional-white noise (\ref{heatEq}). Let  $B^{\frac{1}{2}, \frac{1}{2}}$ be a bifractional Brownian motion with parameters $H=K=\frac{1}{2}$. Consider $(X^{H}_{t})_{t\in [0,T]}$ to be a centered Gaussian process with covariance, for $s,t\in [0,T]$ 
\begin{eqnarray}\label{rx}
R^{X^{H}}(t,s)&=& 2 \frac{1}{2 \sqrt{2\pi}}\alpha _{H} \int_{0} ^{s} (s-a) ^{ 2H-2} \left[ (t+a) ^{\frac{1}{2} } -(t-a) ^{\frac{1}{2}} \right] da\nonumber \\
&=&H \frac{1}{2 \sqrt{2\pi}}\int_{0} ^{s} (s-a) ^{2H-1} \left[ (t+a) ^{-\frac{1}{2} } + (t-a) ^{-\frac{1}{2}}\right] da,
\end{eqnarray}
and let $(Y^{H}_{t}) _{t\in [0,T]}$ be a centered Gaussian process with covariance
\begin{eqnarray}\label{ry}
R^{Y^{H}}(t,s)&=& 2 \frac{1}{2 \sqrt{2\pi}}\alpha _{H}\int_{0} ^{s} (s-a) ^{\frac{1}{2} } \left[ (t+a) ^{2H-2} + (t-a) ^{2H-2} \right] da \nonumber\\
&=&H \frac{1}{ 2\sqrt{2\pi}}\int_{0} ^{s} (s-a) ^{-\frac{1}{2}} \left[ (t+a) ^{2H-1} -(t-a) ^{2H-1} \right] da.
\end{eqnarray}
Suppose that $U, X^{H},$ and $Y^{H}$ are independent. Then for every $x\in \mathbb{R} ^{d}$,
\begin{equation*}
\left( U(t,x) + Y^{H}, t\in [0,T]\right)  \stackrel{\rm Law}{=}\left(  C_{0} B^{\frac{1}{2}, \frac{1}{2}}_{t} + X^{H}_{t}, t\in [0,T]\right),
\end{equation*}where $C_{0}^{2} = \frac{1}{ 2\sqrt{2\pi}}\alpha_{H}C_{d}\beta\left(2H-1,-\frac{d}{2} + 2\right)$.
\end{theorem}
{\bf Proof: }
Let us first justify that $R^{X ^{H}} $ is a covariance function. Clearly, it is symmetric and it can be written, for every $s,t\in [0,T]$, as
\begin{eqnarray*}
\left( \frac{1}{2 \sqrt{2\pi}}\right) ^{-1}R^{X^{H}}(t,s)&=& H\int_{0} ^{s\wedge t} (t\wedge s -a) ^{2H-1} \left[ \left( (t\vee s) +a \right) ^{-\frac{1}{2}} + \left( (t\vee s) -a \right) ^{-\frac{1}{2}}\right] \\
&=& H\int_{0} ^{\infty} 1_{[0,t]} (a) 1_{[0,s]} (a) (t\wedge s -a) ^{2H-1}\left( (t+a) ^{-\frac{1}{2} } \wedge (s+a) ^{-\frac{1}{2}} \right) da \\
&&+ H\int_{0} ^{\infty} 1_{[0,t]} (a) 1_{[0,s]} (a) (t\wedge s -a) ^{2H-1}\left( (t-a) ^{-\frac{1}{2} } \wedge (s-a) ^{-\frac{1}{2}} \right) da
\end{eqnarray*}
and both summands above are positive definite (the same  argument has been used in \cite{BGT}, proof of Theorem 2.1.). Similarly, the function $R^{Y^{H}}$ is a covariance. If $d=1$, we have $C_{d}= 2$ and
\begin{equation*}
R(t,s)= 2 \alpha _{H} \beta (2H-1, \frac{3}{2})\left[ (t+s) ^{\frac{1}{2} } -(t-s) ^{\frac{1}{2} }\right] + R_{1}^{(1)}(t,s)
\end{equation*}
with
\begin{eqnarray*}
\left( \frac{1}{2 \sqrt{2\pi}}\right) ^{-1}R_{1}^{(1)}(t,s)&=&2 \alpha _{H} \int_{0} ^{s} (s-a) ^{ 2H-2} \left[ (t+a) ^{\frac{1}{2} } -(t-a) ^{\frac{1}{2}} \right] da\\
&&-2 \alpha _{H}\int_{0} ^{s} (s-a) ^{\frac{1}{2} } \left[ (t+a) ^{2H-2} + (t-a) ^{2H-2} \right] da\\
&=& H  \int_{0} ^{s} (s-a) ^{ 2H-1}\left[ (t+a) ^{-\frac{1}{2} } +(t-a) ^{-\frac{1}{2}} \right] da\\
&&-2 \alpha _{H}\int_{0} ^{s} (s-a) ^{\frac{1}{2} } \left[ (t+a) ^{2H-2} + (t-a) ^{2H-2} \right] da\\
\end{eqnarray*}
where we used integration by parts in the first integral.
\qed

\noindent In the case $d=3$ we have the following.
\begin{theorem}\label{tt3}
Assume $d=3$. Let $B^{\frac{1}{2}, \frac{1}{2}}$ be a bifractional Brownian motion with $H=K= \frac{1}{2}$ and let $(Z^{H}_{t})_{t\in [0,T]}$ be a centered Gaussian process with covariance $R_{1} ^{(3)} (t,s)$. Then
$$\left( U(t,x)+  C_{0}B^{\frac{1}{2}, \frac{1}{2}}, t\in [0,T] \right) \stackrel{\rm Law}{=} (Z ^{H}_{t}, t\in [0,T]),$$
with $C_{0}$ defined in Theorem \ref{tt2}.

\end{theorem}
{\bf Proof: } We have $C_{3} =-2.$ In this case we can write
\begin{eqnarray*}
R(t,s) + 2 \alpha _{H} \beta (2H-1, \frac{1}{2}) \left[ (t+s) ^{2H-\frac{3}{2}}- (t-s) ^{2H-\frac{3}{2}}\right] =R_{1} ^{(3) } (t,s)
\end{eqnarray*}
with
\begin{eqnarray*}
\left( \frac{1}{2 \sqrt{2\pi}}\right) ^{-1}R_{1} ^{(3)} (t,s)&=& -2\alpha _{H} \int_{0} ^{s} (s-a) ^{2H-2}\left[ (t+a) ^{-\frac{1}{2} } -(t-a) ^{-\frac{1}{2}} \right] da \\
&&+ 2\alpha _{H} \int_{0} ^{s} (s-a) ^{-\frac{1}{2}} \left[ (t+a) ^{2H-2} + (t-a) ^{2H-2} \right].
\end{eqnarray*}

\noindent Note that $R_{1} ^{(3)}$ is a covariance function because it is the sum of two covariance functions. \qed

\begin{remark}\label{onehalf}
Let us understand what happens with the decompositions in Theorems \ref{tt2} and \ref{tt3} when $H $ is  close to $\frac{1}{2}$.  We focus  on the case $d=1$.
The phenomenon is interesting. We first notice that the process $Y^{H}$  vanishes in this case. The covariance of the process $X^{H}$ becomes
$$R^{X^{\frac{1}{2}}}(t,s)=\frac{1}{4\sqrt{2\pi}}\int_{0} ^{s} \left( (t+a) ^{-\frac{1}{2}} -(t-a) ^{-\frac{1}{2}} \right) da=\frac{1}{2 \sqrt{2\pi}}\left( (t+s) ^{\frac{1}{2}}-(t-s) ^{\frac{1}{2}}\right). $$
The constant $C_{0}=\frac{1}{\sqrt{2\pi}}\alpha _{H} \beta (2H-1, \frac{3}{2})$ is not defined for $H=\frac{1}{2}$ because of the presence of $2H-1$ in the argument of the beta function. But the following happens: since $2\alpha _{H} 1_{(0,1) }(u) (1-u) ^{2H-2}$ is an approximation of the unity, it follows that $\alpha _{H} \beta (2H-1, \frac{3}{2})$ converges to $\frac{1}{2}$ when $H$ goes to $\frac{1}{2}$. Therefore $C_{0}$ becomes $\frac{1}{2 \sqrt{2\pi}}$. Therefore we retrieve the result in \cite{Swanson1} and recalled in relation (\ref{covu}). In other words, in the fractional case $H\not= \frac{1}{2}$ the solution ``keeps'' half of the bifractional Brownian motion $B^{\frac{1}{2}, \frac{1}{2}}$ while the other half ``spreads'' into two parts.
\end{remark}

\section{Analysis of the processes $X^{H}$ and $Y^{H}$}

In this paragraph, we will prove some properties of the processes $X^{H}$ and $Y^{H}$ appearing in Theorem \ref{tt2}. It would be nice to find that these processes are more regular than the bifractional Brownian motion $B^{\frac{1}{2}, \frac{1}{2}}$ and  therefore reduce the study of the solution $U$ to the study of $B^{\frac{1}{2}, \frac{1}{2}}$ whose properties are relatively well-known.

\noindent We start with the following remark.
\begin{prop}
The  processes $(X^{H}_{t})_{t\geq 0} , (Y^{H}_{t})_{t\geq 0} $ and $(Z ^{H}_{t})_{t\geq 0} $ from Theorems \ref{tt2} and \ref{tt3} are $2H-\frac{1}{2}$ self-similar.
\end{prop}
{\bf Proof: } Indeed, for every $t\geq s$ and for every $c>0$
\begin{eqnarray*}
\left( \frac{1}{2 \sqrt{2\pi}}\right) ^{-1}R^{X^{H}}(ct, cs) &=& H \int_{0} ^{cs} (cs-a) ^{2H-1} \left[ (ct+a) ^{-\frac{1}{2}}+ (ct-a) ^{-\frac{1}{2}} \right] da \\
&=& H c^{2H-\frac{1}{2}}\left( \frac{1}{2 \sqrt{2\pi}}\right) ^{-1} R^{X^{H}} (t,s)
 \end{eqnarray*}
 where we used the change of variables  $a'= \frac{a}{c}$. The same calculation can be done for $Y^{H}$ and $Z^{H}$. \qed

 \vskip0.3cm

\noindent Let us now estimate the $L^{2}$ norm of the increments of the process $X^{H}$ in order to study the regularity of its paths.

 \begin{prop}For $t>s$ we have
 \begin{equation*}
 \mathbf{E} \left(  \vert X^{H}_{t}-X^{H}_{s} \vert ^{2} \right)  \leq C \vert t-s\vert ^{\frac{1}{2}} \mbox{\  and\  } \mathbf{E} \left(  \vert Y^{H}_{t}-Y^{H}_{s} \vert ^{2} \right)  \leq C \vert t-s\vert ^{\frac{1}{2}}
 \end{equation*}
 with $C$ a strictly positive constant. Consequently, $X^{H}$ and $Y^{H}$ each admit a version that is H\"older continuous of order $\delta \in (0, \frac{1}{4})$.
 \end{prop}
 {\bf Proof: } Using the expression of the covariance of the process $X^{H}$,
\begin{eqnarray*}
\mathbf{E} \vert X^{H}_{t}-X^{H}_{s} \vert ^{2} &=& R^{X^{H}} (t,t)-2R^{X^{H}} (t,s)+ R^{X^{H}} (s,s) \\
&=&  H\int_{0} ^{t} (t-a) ^{2H-1} \left[ (t+a) ^{-\frac{1}{2} } + (t-a) ^{-\frac{1}{2} }\right] da\\
&&+ H\int_{0} ^{a} (s-a) ^{2H-1} \left[ (s+a) ^{-\frac{1}{2} } + (s-a) ^{-\frac{1}{2} }\right] da\\
&&- 2 H\int_{0} ^{s} (s-a) ^{2H-1} \left[ (t+a) ^{-\frac{1}{2} } + (t-a) ^{-\frac{1}{2} }\right] da\\
&:=& T _{1} (t,s) + T_{2}(t,s) + T_{3}(t,s)
\end{eqnarray*}
where
\begin{equation*}
T_{1}(t,s)= H\int_{s}^{t} (t-a) ^{2H-1} \left[ (t+a) ^{-\frac{1}{2} } + (t-a) ^{-\frac{1}{2} }\right] da,
\end{equation*}

\begin{equation*}
T_{2}(t,s) =H\int_{0} ^{s} \left( (t-a) ^{2H-1} - (s-a) ^{2H-1} \right) \left[ (t+a) ^{-\frac{1}{2} } + (t-a) ^{-\frac{1}{2} }\right] da
\end{equation*}
and
\begin{equation*}
T_{3}(t,s)= H\int_{0} ^{s} (s-a)^{2H-1} \left[ (s+a) ^{-\frac{1}{2} } + (s-a) ^{-\frac{1}{2} }-(t+a) ^{-\frac{1}{2} } - (t-a) ^{-\frac{1}{2} }\right]da.
\end{equation*}

\noindent We claim that, for every $t,s \in [0,T]$, $t>s$, we have
\begin{equation*}
C_{1}\vert t-s\vert ^{2H-\frac{1}{2} } \leq T_{1} (t,s) \leq C_{2} \vert t-s \vert  ^{2H-\frac{1}{2} },
\end{equation*}

\begin{equation*}
T_{2} (t,s) \leq C_{3} \vert t-s \vert  ^{2H-\frac{1}{2} }  \mbox{\  and\  }
T_{3} (t,s) \leq C_{4} \vert t-s \vert  ^{\frac{1}{2} }
\end{equation*}
with $C_{1}, C_{2}, C_{3}, C_{4}$ strictly positive constants. Indeed,  since  $(t+a) ^{-\frac{1}{2} } \leq (t-a) ^{-\frac{1}{2}}$ we obtain
\begin{equation*}
T_{1}(t,s) \leq 2H \int_{s} ^{t} (t-a) ^{2H-1} (t-a) ^{-\frac{1}{2}} da = \frac{2H} {2H-\frac{1}{2}}(t-s) ^{2H-\frac{1}{2} }.
\end{equation*}
For the lower bound, since $(t-a) ^{-\frac{1}{2}} \geq (t-s) ^{-\frac{1}{2}}$,
\begin{equation*}
T_{1}(t,s)\geq H(t-s) ^{-\frac{1}{2}}\int_{s}^{t} (t-a) ^{2H-1} da = \frac{H}{2H-1} (t-s) ^{2H-\frac{1}{2} }.
\end{equation*}
Concerning the summands $T_{2}$ and $T_{3}$ we can write
\begin{equation*}
T_{2}(t,s) \leq (t-s) ^{2H-1}H \int_{0} ^{s}(t-a) ^{-\frac{1}{2}} da = \frac{H}{2} (t-s) ^{2H-\frac{1}{2} }.
\end{equation*}
and
\begin{eqnarray*}
T_{3}(t,s) &\leq &H T^{2H-1}\int_{0} ^{s} \left[ (s+a) ^{-\frac{1}{2} } + (s-a) ^{-\frac{1}{2} }-(t+a) ^{-\frac{1}{2} } - (t-a) ^{-\frac{1}{2} }\right]da\\
&=& \frac{H}{2} \left( (2s) ^{\frac{1}{2} } -(t+s) ^{\frac{1}{2} } + (t-s) ^{\frac{1}{2}} \right)\\
&\leq &H  (t-s) ^{\frac{1}{2}}.
\end{eqnarray*}
The calculations are similar for $Y^{H}$. The H\"older continuity of the paths follows by using the Kolmogorov criterium.\qed

\begin{remark}
A  detailed study of the term denoted by $T_{3}$ above shows that its H\"older order cannot be improved.  Therefore the regularity of the paths of $X^{H}$ and $Y^{H}$ is the same as those of the bifractional Brownian motion $B^{\frac{1}{2}, \frac{1}{2}}$. Recall that $B^{\frac{1}{2}, \frac{1}{2}}$ is H\"older continuous of order strictly less than $\frac{1}{4}$.
\end{remark}

\end{document}